\documentclass[12pt]{amsart}
\usepackage{amssymb,epic,eepic}

\def\draft{n}

\theoremstyle{plain}

\newtheorem{theorem}{Theorem}

\newtheorem{proposition}{Proposition}[section]
\newtheorem{lemma}[proposition]{Lemma}

\newtheorem{conjecture}{Conjecture}

\theoremstyle{definition}
\newtheorem{definition}[proposition]{Definition}

\theoremstyle{remark}

\newtheorem{remark}[proposition]{Remark}

\newcommand{\prtag}[1]{\tag*{\llap{$#1$\hskip-\displaywidth}}}

\def\printname#1{
	\if\draft y
		\smash{\makebox[0pt]{\hspace{-0.5in}
			\raisebox{8pt}{\tt\tiny #1}}}
	\fi
}

\newcommand{\mathmode}[1]{$#1$}

\newlength{\standardunitlength}
\def\smily{{
  \setlength{\unitlength}{0.4\standardunitlength}
  \begin{array}{c}  \hspace{-1.7mm} \raisebox{-2pt}{
    \begin{picture}(616,629)(0,-10)
    \thicklines
    \put(308.000,344.500){\arc{375.000}{0.6435}{2.4981}}
    \put(308,307){\ellipse{600}{600}}
    \put(195,382){\blacken\ellipse{76}{76}}
    \put(195,382){\ellipse{76}{76}}
    \put(420,382){\blacken\ellipse{76}{76}}
    \put(420,382){\ellipse{76}{76}}
    \end{picture}
  }
  \hspace{-1.9mm}
  \end{array}
}}

\catcode`\@=11
\long\def\@makecaption#1#2{%
    \vskip 10pt
    \setbox\@tempboxa\hbox{
      \small\sf{\bfcaptionfont #1. }\ignorespaces #2}%
    \ifdim \wd\@tempboxa >\captionwidth {%
        \rightskip=\@captionmargin\leftskip=\@captionmargin
        \unhbox\@tempboxa\par}%
      \else
        \hbox to\hsize{\hfil\box\@tempboxa\hfil}%
    \fi}
\font\bfcaptionfont=cmssbx10 scaled \magstephalf
\newdimen\@captionmargin\@captionmargin=2\parindent
\newdimen\captionwidth\captionwidth=\hsize
\catcode`\@=12

\newcommand{\Diff}{\operatorname{Diff}}

\newcommand{\Span}{\operatorname{span}}
\newcommand{\ad}{\operatorname{ad}}
\newcommand{\gr}{\operatorname{gr}}

\newcommand{\tr}{\operatorname{tr}}

\def\lbl#1{\label{#1}\printname{#1}}

\def\qed{{\hfill\text{$\Box$}}}

\newenvironment{myitemize}{
	\begin{list}{$\bullet$}{\setlength{\leftmargin}{16pt}
        \setlength{\labelwidth}{12pt}
        \setlength{\labelsep}{4pt}}
}{
	\end{list}
}

\newlength{\globalparindent}
\newcommand{\udot}{{\mathaccent\cdot\cup}}

\advance \topmargin by -\headheight
\advance \topmargin by -\headsep
\evensidemargin \oddsidemargin
\newcommand{\Arhus}{\AA rhus}

\def\calA{{\mathcal A}}
\def\calB{{\mathcal B}}

\def\calV{{\mathcal V}}
\def\frakg{{\mathfrak g}}
\def\frakh{{\mathfrak h}}

\def\KF{{{\mathcal K}^F}}
\def\ZF{Z}
\def\RTg{{RT_\frakg}}
\def\A{{\calA'}}
\def\Tg{{{\mathcal T}_\frakg^\hbar}}
\def\Ug{{U(\frakg)^\frakg[[\hbar]]}}
\def\Ph{{P(\frakh^\star)^W[[\hbar]]}}
\def\ch{{\chi}}
\def\b{{\beta_\frakg}}
\def\psig{{\psi_\frakg}}
\def\Bt{{\calB'_\times}}
\def\Bu{{\calB'_\udot}}
\def\St{{S(\frakg)^\frakg_\times[[\hbar]]}}
\def\Su{{S(\frakg)^\frakg_\udot[[\hbar]]}}
\def\Pg{{P(\frakg^\star)^\frakg}[[\hbar]]}
\def\Oh{{\hat{\Omega}}}
\def\Dj{{D(j_\frakg^{1/2})}}
\def\ig{{\iota_\frakg}}
\def\Sg{{S_\frakg}}
\def\Ca{{\setlength{\unitlength}{0.6\standardunitlength}
	\begin{array}{c}  \hspace{-1.7mm}
        	\raisebox{-2pt}{%
\begingroup\makeatletter\ifx\SetFigFont\undefined
\def\x#1#2#3#4#5#6#7\relax{\def\x{#1#2#3#4#5#6}}%
\expandafter\x\fmtname xxxxxx\relax \def\y{splain}%
\ifx\x\y   
\gdef\SetFigFont#1#2#3{%
  \ifnum #1<17\tiny\else \ifnum #1<20\small\else
  \ifnum #1<24\normalsize\else \ifnum #1<29\large\else
  \ifnum #1<34\Large\else \ifnum #1<41\LARGE\else
     \huge\fi\fi\fi\fi\fi\fi
  \csname #3\endcsname}%
\else
\gdef\SetFigFont#1#2#3{\begingroup
  \count@#1\relax \ifnum 25<\count@\count@25\fi
  \def\x{\endgroup\@setsize\SetFigFont{#2pt}}%
  \expandafter\x
    \csname \romannumeral\the\count@ pt\expandafter\endcsname
    \csname @\romannumeral\the\count@ pt\endcsname
  \csname #3\endcsname}%
\fi
\fi\endgroup
{\renewcommand{\dashlinestretch}{30}
\begin{picture}(474,414)(0,-10)
\path(12,387)(12,12)(462,12)(12,387)
\put(87,87){\makebox(0,0)[lb]{\smash{{\mathmode{\scriptstyle 1}}}}}
\end{picture}
}
 }
        	\hspace{-1.9mm}
	\end{array}
}}
\def\Cb{{\setlength{\unitlength}{0.5\standardunitlength}
	\begin{array}{c}  \hspace{-1.7mm}
        	\raisebox{-2pt}{%
\begingroup\makeatletter\ifx\SetFigFont\undefined
\def\x#1#2#3#4#5#6#7\relax{\def\x{#1#2#3#4#5#6}}%
\expandafter\x\fmtname xxxxxx\relax \def\y{splain}%
\ifx\x\y   
\gdef\SetFigFont#1#2#3{%
  \ifnum #1<17\tiny\else \ifnum #1<20\small\else
  \ifnum #1<24\normalsize\else \ifnum #1<29\large\else
  \ifnum #1<34\Large\else \ifnum #1<41\LARGE\else
     \huge\fi\fi\fi\fi\fi\fi
  \csname #3\endcsname}%
\else
\gdef\SetFigFont#1#2#3{\begingroup
  \count@#1\relax \ifnum 25<\count@\count@25\fi
  \def\x{\endgroup\@setsize\SetFigFont{#2pt}}%
  \expandafter\x
    \csname \romannumeral\the\count@ pt\expandafter\endcsname
    \csname @\romannumeral\the\count@ pt\endcsname
  \csname #3\endcsname}%
\fi
\fi\endgroup
{\renewcommand{\dashlinestretch}{30}
\begin{picture}(549,414)(0,-10)
\path(12,387)(537,387)(537,12)
	(12,12)(12,387)
\put(192,132){\makebox(0,0)[lb]{\smash{{\mathmode{\scriptstyle 2}}}}}
\end{picture}
}
 }
        	\hspace{-1.9mm}
	\end{array}
}}
\def\Cc{{\setlength{\unitlength}{0.4\standardunitlength}
	\begin{array}{c}  \hspace{-1.7mm}
        	\raisebox{-2pt}{%
\begingroup\makeatletter\ifx\SetFigFont\undefined
\def\x#1#2#3#4#5#6#7\relax{\def\x{#1#2#3#4#5#6}}%
\expandafter\x\fmtname xxxxxx\relax \def\y{splain}%
\ifx\x\y   
\gdef\SetFigFont#1#2#3{%
  \ifnum #1<17\tiny\else \ifnum #1<20\small\else
  \ifnum #1<24\normalsize\else \ifnum #1<29\large\else
  \ifnum #1<34\Large\else \ifnum #1<41\LARGE\else
     \huge\fi\fi\fi\fi\fi\fi
  \csname #3\endcsname}%
\else
\gdef\SetFigFont#1#2#3{\begingroup
  \count@#1\relax \ifnum 25<\count@\count@25\fi
  \def\x{\endgroup\@setsize\SetFigFont{#2pt}}%
  \expandafter\x
    \csname \romannumeral\the\count@ pt\expandafter\endcsname
    \csname @\romannumeral\the\count@ pt\endcsname
  \csname #3\endcsname}%
\fi
\fi\endgroup
{\renewcommand{\dashlinestretch}{30}
\begin{picture}(849,789)(0,-10)
\path(12,762)(12,387)(312,12)
	(837,12)(612,762)(12,762)
\put(312,312){\makebox(0,0)[lb]{\smash{{\mathmode{\scriptstyle 3}}}}}
\end{picture}
}
 }
        	\hspace{-1.9mm}
	\end{array}
}}
\def\Cd{{\setlength{\unitlength}{0.5\standardunitlength}
	\begin{array}{c}  \hspace{-1.7mm}
        	\raisebox{-2pt}{%
\begingroup\makeatletter\ifx\SetFigFont\undefined
\def\x#1#2#3#4#5#6#7\relax{\def\x{#1#2#3#4#5#6}}%
\expandafter\x\fmtname xxxxxx\relax \def\y{splain}%
\ifx\x\y   
\gdef\SetFigFont#1#2#3{%
  \ifnum #1<17\tiny\else \ifnum #1<20\small\else
  \ifnum #1<24\normalsize\else \ifnum #1<29\large\else
  \ifnum #1<34\Large\else \ifnum #1<41\LARGE\else
     \huge\fi\fi\fi\fi\fi\fi
  \csname #3\endcsname}%
\else
\gdef\SetFigFont#1#2#3{\begingroup
  \count@#1\relax \ifnum 25<\count@\count@25\fi
  \def\x{\endgroup\@setsize\SetFigFont{#2pt}}%
  \expandafter\x
    \csname \romannumeral\the\count@ pt\expandafter\endcsname
    \csname @\romannumeral\the\count@ pt\endcsname
  \csname #3\endcsname}%
\fi
\fi\endgroup
{\renewcommand{\dashlinestretch}{30}
\begin{picture}(849,414)(0,-10)
\path(12,387)(537,387)(837,12)
	(312,12)(12,387)
\put(342,132){\makebox(0,0)[lb]{\smash{{\mathmode{\scriptstyle 4}}}}}
\end{picture}
}
 }
        	\hspace{-1.9mm}
	\end{array}
}}
\def\Ce{{\setlength{\unitlength}{0.5\standardunitlength}
	\begin{array}{c}  \hspace{-1.7mm}
        	\raisebox{-2pt}{%
\begingroup\makeatletter\ifx\SetFigFont\undefined
\def\x#1#2#3#4#5#6#7\relax{\def\x{#1#2#3#4#5#6}}%
\expandafter\x\fmtname xxxxxx\relax \def\y{splain}%
\ifx\x\y   
\gdef\SetFigFont#1#2#3{%
  \ifnum #1<17\tiny\else \ifnum #1<20\small\else
  \ifnum #1<24\normalsize\else \ifnum #1<29\large\else
  \ifnum #1<34\Large\else \ifnum #1<41\LARGE\else
     \huge\fi\fi\fi\fi\fi\fi
  \csname #3\endcsname}%
\else
\gdef\SetFigFont#1#2#3{\begingroup
  \count@#1\relax \ifnum 25<\count@\count@25\fi
  \def\x{\endgroup\@setsize\SetFigFont{#2pt}}%
  \expandafter\x
    \csname \romannumeral\the\count@ pt\expandafter\endcsname
    \csname @\romannumeral\the\count@ pt\endcsname
  \csname #3\endcsname}%
\fi
\fi\endgroup
{\renewcommand{\dashlinestretch}{30}
\begin{picture}(849,414)(0,-10)
\path(12,387)(537,387)(837,12)
	(312,12)(12,387)
\put(342,132){\makebox(0,0)[lb]{\smash{{\mathmode{\scriptstyle 5}}}}}
\end{picture}
}
 }
        	\hspace{-1.9mm}
	\end{array}
}}

\begin{document}

\title{Wheels, Wheeling, and the Kontsevich Integral of the Unknot}

\author[Bar-Natan]{Dror~Bar-Natan}
\address{Institute of Mathematics\\
        The Hebrew University\\
        Giv'at-Ram, Jerusalem 91904\\
        Israel}
\email{drorbn@math.huji.ac.il}

\author[Garoufalidis]{Stavros~Garoufalidis}
\address{Department of Mathematics\\
        Harvard University\\
        Cambridge MA 02138\\
        USA}
\email{stavros@math.harvard.edu}

\author[Rozansky]{Lev~Rozansky}
\address{Department of Mathematics, Statistics, and Computer Science\\
        University of Illinois at Chicago\\
        Chicago IL 60607-7045\\
        USA}
\email{rozansky@math.uic.edu}

\author[Thurston]{Dylan~P.~Thurston}
\address{Department of Mathematics\\
        University of California at Berkeley\\
        Berkeley CA 94720-3840\\
        USA}
\email{dpt@math.berkeley.edu}

\thanks{This preprint is available electronically at
  {\tt http://www.ma.huji.ac.il/\~{}drorbn}, at \newline
  {\tt http://jacobi.math.brown.edu/\~{}stavrosg}, and at {\tt
    http://xxx.lanl.gov/abs/q-alg/9703025}.
}

\dedicatory{This is a preprint. Your comments are welcome.}
\date{This edition: Apr.~26,~1998; \ \ First edition: Mar.~13, 1997.}

\begin{abstract}
We conjecture an exact formula for the Kontsevich integral of the
unknot, and also conjecture a formula (also conjectured independently by
Deligne~\cite{Deligne:Letter}) for the relation between the two natural
products on the space of Chinese characters. The two formulas use the
related notions of ``Wheels'' and ``Wheeling''.  We prove these formulas
`on the level of Lie algebras' using standard techniques from the theory
of Vassiliev invariants and the theory of Lie algebras.
\end{abstract}

\maketitle

\tableofcontents

\section{Introduction}

\subsection{The conjectures}
Let us start with the statements of our conjectures; the rest of
the paper is concerned with motivating and justifying them. We
assume some familiarity with the theory of Vassiliev invariants. See
e.g.~\cite{Bar-Natan:Vassiliev, Birman:Bulletin, BirmanLin:Vassiliev,
Goussarov:New, Goussarov:nEquivalence, Kontsevich:Vassiliev,
Vassiliev:CohKnot, Vassiliev:Book} and~\cite{Bar-Natan:VasBib}.

Very briefly, recall that any complex-valued knot invariant $V$ can be
extended to an invariant of knots with double points ({\em singular
knots}) via the formula $V(\!\setlength{\unitlength}{0.5\standardunitlength}
	\begin{array}{c}  \hspace{-1.7mm}
        	\raisebox{-2pt}{%
\begingroup\makeatletter\ifx\SetFigFont\undefined
\def\x#1#2#3#4#5#6#7\relax{\def\x{#1#2#3#4#5#6}}%
\expandafter\x\fmtname xxxxxx\relax \def\y{splain}%
\ifx\x\y   
\gdef\SetFigFont#1#2#3{%
  \ifnum #1<17\tiny\else \ifnum #1<20\small\else
  \ifnum #1<24\normalsize\else \ifnum #1<29\large\else
  \ifnum #1<34\Large\else \ifnum #1<41\LARGE\else
     \huge\fi\fi\fi\fi\fi\fi
  \csname #3\endcsname}%
\else
\gdef\SetFigFont#1#2#3{\begingroup
  \count@#1\relax \ifnum 25<\count@\count@25\fi
  \def\x{\endgroup\@setsize\SetFigFont{#2pt}}%
  \expandafter\x
    \csname \romannumeral\the\count@ pt\expandafter\endcsname
    \csname @\romannumeral\the\count@ pt\endcsname
  \csname #3\endcsname}%
\fi
\fi\endgroup
{\renewcommand{\dashlinestretch}{30}
\begin{picture}(324,339)(0,-10)
\put(162,162){\blacken\ellipse{80}{80}}
\put(162,162){\ellipse{80}{80}}
\path(12,12)(312,312)
\path(248.360,205.934)(312.000,312.000)(205.934,248.360)
\path(312,12)(12,312)
\path(118.066,248.360)(12.000,312.000)(75.640,205.934)
\end{picture}
}
 }
        	\hspace{-1.9mm}
	\end{array}
\!\!) =
V(\!\setlength{\unitlength}{0.5\standardunitlength}
	\begin{array}{c}  \hspace{-1.7mm}
        	\raisebox{-2pt}{%
\begingroup\makeatletter\ifx\SetFigFont\undefined
\def\x#1#2#3#4#5#6#7\relax{\def\x{#1#2#3#4#5#6}}%
\expandafter\x\fmtname xxxxxx\relax \def\y{splain}%
\ifx\x\y   
\gdef\SetFigFont#1#2#3{%
  \ifnum #1<17\tiny\else \ifnum #1<20\small\else
  \ifnum #1<24\normalsize\else \ifnum #1<29\large\else
  \ifnum #1<34\Large\else \ifnum #1<41\LARGE\else
     \huge\fi\fi\fi\fi\fi\fi
  \csname #3\endcsname}%
\else
\gdef\SetFigFont#1#2#3{\begingroup
  \count@#1\relax \ifnum 25<\count@\count@25\fi
  \def\x{\endgroup\@setsize\SetFigFont{#2pt}}%
  \expandafter\x
    \csname \romannumeral\the\count@ pt\expandafter\endcsname
    \csname @\romannumeral\the\count@ pt\endcsname
  \csname #3\endcsname}%
\fi
\fi\endgroup
{\renewcommand{\dashlinestretch}{30}
\begin{picture}(324,339)(0,-10)
\path(312,12)(237,87)
\path(87,237)(12,312)
\path(118.066,248.360)(12.000,312.000)(75.640,205.934)
\path(12,12)(312,312)
\path(248.360,205.934)(312.000,312.000)(205.934,248.360)
\end{picture}
}
 }
        	\hspace{-1.9mm}
	\end{array}
\!\!) -
V(\!\setlength{\unitlength}{0.5\standardunitlength}
	\begin{array}{c}  \hspace{-1.7mm}
        	\raisebox{-2pt}{%
\begingroup\makeatletter\ifx\SetFigFont\undefined
\def\x#1#2#3#4#5#6#7\relax{\def\x{#1#2#3#4#5#6}}%
\expandafter\x\fmtname xxxxxx\relax \def\y{splain}%
\ifx\x\y   
\gdef\SetFigFont#1#2#3{%
  \ifnum #1<17\tiny\else \ifnum #1<20\small\else
  \ifnum #1<24\normalsize\else \ifnum #1<29\large\else
  \ifnum #1<34\Large\else \ifnum #1<41\LARGE\else
     \huge\fi\fi\fi\fi\fi\fi
  \csname #3\endcsname}%
\else
\gdef\SetFigFont#1#2#3{\begingroup
  \count@#1\relax \ifnum 25<\count@\count@25\fi
  \def\x{\endgroup\@setsize\SetFigFont{#2pt}}%
  \expandafter\x
    \csname \romannumeral\the\count@ pt\expandafter\endcsname
    \csname @\romannumeral\the\count@ pt\endcsname
  \csname #3\endcsname}%
\fi
\fi\endgroup
{\renewcommand{\dashlinestretch}{30}
\begin{picture}(324,339)(0,-10)
\path(12,12)(87,87)
\path(237,237)(312,312)
\path(248.360,205.934)(312.000,312.000)(205.934,248.360)
\path(312,12)(12,312)
\path(118.066,248.360)(12.000,312.000)(75.640,205.934)
\end{picture}
}
 }
        	\hspace{-1.9mm}
	\end{array}
\!\!)$.  An invariant of knots
(or framed knots) is called a {\em Vassiliev invariant}, or a {\em
finite type invariant of type $m$}, if its extension to singular knots
vanishes whenever evaluated on a singular knot that has more than $m$
double points.  Vassiliev invariants are in some senses analogues to
polynomials (on the space of all knots), and one may hope that they
separate knots. While this is an open problem and the precise power of the
Vassiliev theory is yet unknown, it is known (see~\cite{Vogel:Structures})
that Vassiliev invariants are strictly stronger than the
Reshetikhin-Turaev invariants (\cite{ReshetikhinTuraev:Ribbon}), and in
particular they are strictly stronger than the Alexander-Conway, Jones,
HOMFLY, and Kauffman invariants. Hence one is interested in a detailed
understanding of the theory of Vassiliev invariants.

The set $\calV$ of all Vassiliev invariants of framed knots is a linear
space, filtered by the ``type'' of an invariant. The fundamental theorem
of Vassiliev invariants, due to Kontsevich~\cite{Kontsevich:Vassiliev},
says that the associated graded space $\gr\calV$ of $\calV$ can be
identified with the graded dual $\calA^\star$ of a certain completed
graded space $\calA$ \lbl{Adefinition} of formal linear combinations of
certain diagrams, modulo certain linear relations. The ``diagrams'' in
$\calA$ are connected graphs made of a single distinguished directed line
(the ``skeleton''), some number of undirected ``internal edges'', some
number of trivalent ``external vertices'' in which an internal edge ends
on the skeleton, and some number of trivalent ``internal vertices''
in which three internal edges meet. It is further assumed that the
internal vertices are ``oriented''; that for each internal vertices
one of the two possible cyclic orderings of the edges emanating
from it is specified. An example of a diagram in $\calA$ is in
figure~\ref{BasicDefinitions}. The linear relations in the definition
of $\calA$ are the well-known $AS$, $IHX$, and $STU$ relations, also
shown in figure~\ref{BasicDefinitions}. The space $\calA$ is graded by
half the total number of trivalent vertices in a given diagram.

\begin{figure}[htpb]
\[ \printname{BasicDefinitions}
	\setlength{\unitlength}{0.85\standardunitlength}
	\begin{array}{c}  \hspace{-1.7mm}
        	\raisebox{-8pt}{\begingroup\makeatletter\ifx\SetFigFont\undefined
\def\x#1#2#3#4#5#6#7\relax{\def\x{#1#2#3#4#5#6}}%
\expandafter\x\fmtname xxxxxx\relax \def\y{splain}%
\ifx\x\y   
\gdef\SetFigFont#1#2#3{%
  \ifnum #1<17\tiny\else \ifnum #1<20\small\else
  \ifnum #1<24\normalsize\else \ifnum #1<29\large\else
  \ifnum #1<34\Large\else \ifnum #1<41\LARGE\else
     \huge\fi\fi\fi\fi\fi\fi
  \csname #3\endcsname}%
\else
\gdef\SetFigFont#1#2#3{\begingroup
  \count@#1\relax \ifnum 25<\count@\count@25\fi
  \def\x{\endgroup\@setsize\SetFigFont{#2pt}}%
  \expandafter\x
    \csname \romannumeral\the\count@ pt\expandafter\endcsname
    \csname @\romannumeral\the\count@ pt\endcsname
  \csname #3\endcsname}%
\fi
\fi\endgroup
{\renewcommand{\dashlinestretch}{30}
\begin{picture}(6024,1257)(0,-10)
\put(912.000,462.000){\arc{900.000}{3.1416}{6.2832}}
\put(462.000,462.000){\arc{600.000}{3.1416}{6.2832}}
\put(3312,762){\ellipse{300}{300}}
\path(12,462)(1512,462)
\path(1392.000,432.000)(1512.000,462.000)(1392.000,492.000)
\path(912,912)(1062,462)
\path(2412,912)(2862,912)(2862,612)
	(2412,612)(2412,912)
\path(2262,1062)(2412,912)
\path(2862,912)(3012,1062)
\path(2862,612)(3012,462)
\path(2412,612)(2262,462)
\path(2637,912)(2637,612)
\path(3312,1062)(3312,912)
\path(3312,612)(3312,462)
\path(4812,762)(5112,762)
\path(4962,762)(4962,462)
\path(4812,462)(5112,462)
\path(5262,762)(5262,462)
\path(5262,612)(5562,612)
\path(5562,762)(5562,462)
\path(5712,762)(6012,462)
\path(6012,762)(5712,462)
\path(5787,537)(5937,537)
\path(4812,12)(5112,12)
\path(4992.000,-18.000)(5112.000,12.000)(4992.000,42.000)
\path(5262,12)(5562,12)
\path(5442.000,-18.000)(5562.000,12.000)(5442.000,42.000)
\path(5712,12)(6012,12)
\path(5892.000,-18.000)(6012.000,12.000)(5892.000,42.000)
\path(4812,312)(4962,162)(5112,312)
\path(4962,162)(4962,12)
\path(5262,312)(5337,12)
\path(5487,312)(5412,12)
\path(5712,312)(5937,12)
\path(6012,312)(5787,12)
\path(4812,1212)(4962,1062)(5112,1212)
\path(4962,1062)(4962,912)
\path(5412,912)(5412,987)
\path(5562,1212)(5561,1212)(5558,1210)
	(5551,1207)(5539,1201)(5522,1193)
	(5503,1184)(5484,1175)(5465,1166)
	(5448,1158)(5434,1150)(5422,1143)
	(5412,1137)(5401,1129)(5392,1122)
	(5382,1114)(5373,1107)(5364,1099)
	(5356,1092)(5348,1084)(5342,1077)
	(5338,1069)(5337,1062)(5338,1054)
	(5343,1046)(5349,1037)(5358,1027)
	(5367,1017)(5377,1008)(5387,999)
	(5396,993)(5404,988)(5412,987)
	(5420,988)(5428,993)(5437,999)
	(5447,1008)(5457,1017)(5466,1027)
	(5475,1037)(5481,1046)(5486,1054)
	(5487,1062)(5486,1070)(5482,1077)
	(5476,1085)(5468,1092)(5460,1100)
	(5451,1107)(5442,1115)(5432,1122)
	(5423,1130)(5412,1137)(5402,1143)
	(5390,1150)(5376,1158)(5359,1166)
	(5340,1175)(5321,1185)(5302,1193)
	(5285,1201)(5273,1207)(5266,1210)
	(5263,1212)(5262,1212)
\put(312,237){\makebox(0,0)[lb]{\smash{{{\rm\scriptsize degree=3}}}}}
\put(2412,237){\makebox(0,0)[lb]{\smash{{{\rm\scriptsize degree=7}}}}}
\put(4212,987){\makebox(0,0)[lb]{\smash{{{\rm\scriptsize AS:}}}}}
\put(4212,537){\makebox(0,0)[lb]{\smash{{{\rm\scriptsize IHX:}}}}}
\put(4212,87){\makebox(0,0)[lb]{\smash{{{\rm\scriptsize STU:}}}}}
\put(5112,162){\makebox(0,0)[lb]{\smash{{\mathmode{\scriptstyle =}}}}}
\put(5112,612){\makebox(0,0)[lb]{\smash{{\mathmode{\scriptstyle =}}}}}
\put(5637,1062){\makebox(0,0)[lb]{\smash{{\mathmode{\scriptstyle =0}}}}}
\put(5112,1062){\makebox(0,0)[lb]{\smash{{\mathmode{\scriptstyle +}}}}}
\put(5562,612){\makebox(0,0)[lb]{\smash{{\mathmode{\scriptstyle -}}}}}
\put(5562,162){\makebox(0,0)[lb]{\smash{{\mathmode{\scriptstyle -}}}}}
\end{picture}
} }
        	\hspace{-1.9mm}
	\end{array}
 \]
\caption{
  A diagram in $\calA$, a diagram in $\calB$ (a Chinese character),
  and the $AS$, $IHX$, and $STU$ relations. All internal vertices shown
  are oriented counterclockwise.
}
\lbl{BasicDefinitions}
\end{figure}
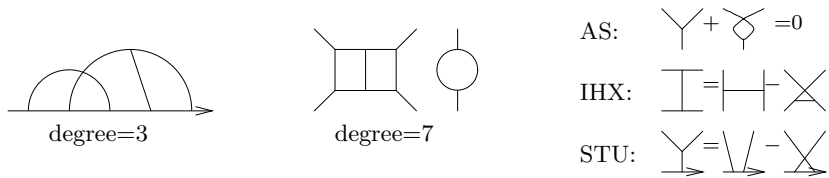

The most difficult part of the currently known proofs of the
isomorphism $\gr\calV\cong\calA^\star$ is the construction of a
``universal Vassiliev invariant''; an $\calA$-valued framed-knot
invariant that satisfies a certain universality property. Such a
``universal Vassiliev invariant'' is not unique; the set of universal
Vassiliev invariants is in a bijective correspondence with the set of all
filtration-respecting maps $\calV\to\gr\calV$ that induce the identity
map $\gr\calV\to\gr\calV$. But it is a noteworthy and not terribly well
understood fact that all known constructions of a universal Vassiliev
invariant are either known to give the same answer or are conjectured to
give the same answer as the original ``framed Kontsevich integral'' $\ZF$
(see Section~\ref{TheEdges}). Furthermore, the Kontsevich integral is
well behaved in several senses, as shown in~\cite{Bar-Natan:Vassiliev,
Bar-NatanGaroufalidis:MMR, Kassel:QuantumGroups, Kontsevich:Vassiliev,
LeMurakami2Ohtsuki:3Manifold, LeMurakami:Universal, LeMurakami:Parallel}.

Thus it seems that $\ZF$ is a canonical and not an accidental object. It
is therefore surprising how little we know about it. While there are
several formulas for computing $\ZF$, they are all of limited use beyond
the first few degrees. Presently, we do not know how to compute $\ZF$
for {\em any} knot; not even the unknot!

Our first conjecture is about the value of the Kontsevich integral of the
unknot. We conjecture a completely explicit formula, written in terms
of an alternative realization of the space $\calA$, the space $\calB$
of ``Chinese characters'' (see~\cite{Bar-Natan:Vassiliev}). The space
$\calB$ is also a completed graded space of formal linear combinations
of diagrams modulo linear relations: the diagrams are the so-called
Chinese characters, which are the same as the diagrams in $\calA$
except that a skeleton is not present, and instead a certain number
of univalent vertices are allowed (the connectivity requirement is
dropped, but one insists that every connected component of a Chinese
character would have at least one univalent vertex).  An example of a
Chinese character is in figure~\ref{BasicDefinitions}. The relations
are the $AS$ and $IHX$ relations that appear in the same figure (but
not the $STU$ relation, which involves the skeleton). The degree of a
Chinese character is half the total number of its vertices. There is a
natural isomorphism $\chi:\calB\to\calA$ \lbl{chidefinition} which maps
every Chinese character to the average of all possible ways of placing
its univalent vertices along a skeleton line. In a sense that we will
recall below, the fact that $\chi$ is an isomorphism is an analog of
the Poincare-Birkhoff-Witt (PBW) theorem. We note that the inverse map
$\sigma$ of $\chi$ is more difficult to construct and manipulate.

\begin{conjecture} \lbl{UnknotConjecture} (Wheels)
The framed Kontsevich integral of the unknot,
$\ZF(\bigcirc)$, expressed in terms of Chinese characters, is equal to
\begin{equation} \lbl{OmegaDef}
  \Omega=\exp_\udot \sum_{n=1}^\infty b_{2n}\omega_{2n}.
\end{equation}
\end{conjecture}

The notation in~\eqref{OmegaDef} means:
\begin{myitemize}
\item The `modified Bernoulli numbers' $b_{2n}$ are defined by the power
series expansion
\begin{equation} \lbl{MBNDefinition}
  \sum_{n=0}^\infty b_{2n}x^{2n} = \frac{1}{2}\log\frac{\sinh x/2}{x/2}.
\end{equation}
These numbers are related to the usual Bernoulli numbers $B_{2n}$ and
to the values of the Riemann $\zeta$-function on the even integers via
(see e.g.~\cite[Section~12.12]{Apostol:AnalyticNumberTheory})
\[
  b_{2n}
  =\frac{B_{2n}}{4n(2n)!}
  =\frac{(-1)^{n+1}}{2n(2\pi)^{2n}}\zeta(2n).
\]
The first three modified Bernoulli numbers are $b_2=1/48$, $b_4=-1/5760$,
and $b_6=1/362880$.
\item The `$2n$-wheel' $\omega_{2n}$ is the degree $2n$ Chinese character
made of a $2n$-gon with $2n$ legs:
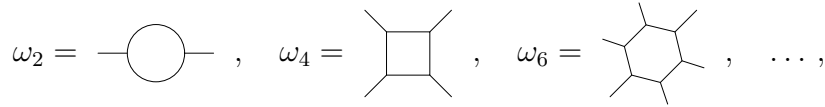
\[
  \omega_2=\printname{2wheel}
	\setlength{\unitlength}{0.6\standardunitlength}
	\begin{array}{c}  \hspace{-1.7mm}
        	\raisebox{-8pt}{%
\begingroup\makeatletter\ifx\SetFigFont\undefined
\def\x#1#2#3#4#5#6#7\relax{\def\x{#1#2#3#4#5#6}}%
\expandafter\x\fmtname xxxxxx\relax \def\y{splain}%
\ifx\x\y   
\gdef\SetFigFont#1#2#3{%
  \ifnum #1<17\tiny\else \ifnum #1<20\small\else
  \ifnum #1<24\normalsize\else \ifnum #1<29\large\else
  \ifnum #1<34\Large\else \ifnum #1<41\LARGE\else
     \huge\fi\fi\fi\fi\fi\fi
  \csname #3\endcsname}%
\else
\gdef\SetFigFont#1#2#3{\begingroup
  \count@#1\relax \ifnum 25<\count@\count@25\fi
  \def\x{\endgroup\@setsize\SetFigFont{#2pt}}%
  \expandafter\x
    \csname \romannumeral\the\count@ pt\expandafter\endcsname
    \csname @\romannumeral\the\count@ pt\endcsname
  \csname #3\endcsname}%
\fi
\fi\endgroup
{\renewcommand{\dashlinestretch}{30}
\begin{picture}(1224,629)(0,-10)
\put(612,307){\ellipse{600}{600}}
\path(12,307)(312,307)
\path(912,307)(1212,307)
\end{picture}
}
 }
        	\hspace{-1.9mm}
	\end{array}
,\quad
  \omega_4=\printname{4wheel}
	\setlength{\unitlength}{0.6\standardunitlength}
	\begin{array}{c}  \hspace{-1.7mm}
        	\raisebox{-8pt}{%
\begingroup\makeatletter\ifx\SetFigFont\undefined
\def\x#1#2#3#4#5#6#7\relax{\def\x{#1#2#3#4#5#6}}%
\expandafter\x\fmtname xxxxxx\relax \def\y{splain}%
\ifx\x\y   
\gdef\SetFigFont#1#2#3{%
  \ifnum #1<17\tiny\else \ifnum #1<20\small\else
  \ifnum #1<24\normalsize\else \ifnum #1<29\large\else
  \ifnum #1<34\Large\else \ifnum #1<41\LARGE\else
     \huge\fi\fi\fi\fi\fi\fi
  \csname #3\endcsname}%
\else
\gdef\SetFigFont#1#2#3{\begingroup
  \count@#1\relax \ifnum 25<\count@\count@25\fi
  \def\x{\endgroup\@setsize\SetFigFont{#2pt}}%
  \expandafter\x
    \csname \romannumeral\the\count@ pt\expandafter\endcsname
    \csname @\romannumeral\the\count@ pt\endcsname
  \csname #3\endcsname}%
\fi
\fi\endgroup
{\renewcommand{\dashlinestretch}{30}
\begin{picture}(924,939)(0,-10)
\path(237,687)(687,687)(687,237)
	(237,237)(237,687)
\path(687,687)(912,912)
\path(687,237)(912,12)
\path(237,237)(12,12)
\path(237,687)(12,912)
\end{picture}
}
 }
        	\hspace{-1.9mm}
	\end{array}
,\quad
  \omega_6=\printname{6wheel}
	\setlength{\unitlength}{0.6\standardunitlength}
	\begin{array}{c}  \hspace{-1.7mm}
        	\raisebox{-8pt}{%
\begingroup\makeatletter\ifx\SetFigFont\undefined
\def\x#1#2#3#4#5#6#7\relax{\def\x{#1#2#3#4#5#6}}%
\expandafter\x\fmtname xxxxxx\relax \def\y{splain}%
\ifx\x\y   
\gdef\SetFigFont#1#2#3{%
  \ifnum #1<17\tiny\else \ifnum #1<20\small\else
  \ifnum #1<24\normalsize\else \ifnum #1<29\large\else
  \ifnum #1<34\Large\else \ifnum #1<41\LARGE\else
     \huge\fi\fi\fi\fi\fi\fi
  \csname #3\endcsname}%
\else
\gdef\SetFigFont#1#2#3{\begingroup
  \count@#1\relax \ifnum 25<\count@\count@25\fi
  \def\x{\endgroup\@setsize\SetFigFont{#2pt}}%
  \expandafter\x
    \csname \romannumeral\the\count@ pt\expandafter\endcsname
    \csname @\romannumeral\the\count@ pt\endcsname
  \csname #3\endcsname}%
\fi
\fi\endgroup
{\renewcommand{\dashlinestretch}{30}
\begin{picture}(1074,1089)(0,-10)
\path(762,762)(837,462)(612,237)
	(312,312)(237,612)(462,837)
	(762,762)(762,762)
\path(762,762)(987,987)
\path(837,462)(1062,387)
\path(612,237)(687,12)
\path(312,312)(87,87)
\path(237,612)(12,687)
\path(462,837)(387,1062)
\end{picture}
}
 }
        	\hspace{-1.9mm}
	\end{array}
,\quad\ldots,
\]
(with all vertices oriented counterclockwise).\footnote{
  Wheels have appeared in several noteworthy places before:
  \cite{Chmutov:CombinatorialMelvinMorton, ChmutovVarchenko:sl2,
  KrickerSpenceAitchison:Cabling, Vaintrob:Primitive}. Similar but slightly
  different objects appear in Ng's beautiful work on ribbon
  knots~\cite{Ng:Ribbon}.
}

\item $\exp_\udot$ means `exponential in the disjoint union sense'; that is,
  it is the formal-sum exponential of a linear combination of Chinese
  characters, with the product being the disjoint union product.
\end{myitemize}

Let us explain why we believe the Wheels Conjecture
(Conjecture~\ref{UnknotConjecture}). Recall (\cite{Bar-Natan:Vassiliev}) that
there is a parallelism between the space $\calA$ (and various variations
thereof) and a certain part of the theory of Lie algebras. Specifically,
given a metrized Lie algebra $\frakg$, there exists a commutative square
(a refined version is in Theorem~\ref{CommutativityTheorem} below)
\[ {
  \def\A{{\mathcal A}}
  \def\B{{\mathcal B}}
  \def\Tg{{{\mathcal T}_\frakg}}
  \def\U{{{U}^\frakg(\frakg})}
  \def\S{{{S}^\frakg(\frakg})}
  \def\Udef{{\left(\parbox{2.4in}{\small the $\frakg$-invariant part
    of the completed universal enveloping algebra of $\frakg$}
  \right)}}
  \def\Sdef{{\left(\parbox{2.1in}{\small the $\frakg$-invariant part
    of the completed symmetric algebra of $\frakg$}
  \right)}}
  \printname{CommutativeSquare}
	\setlength{\unitlength}{0.8\standardunitlength}
	\begin{array}{c}  \hspace{-1.7mm}
        	\raisebox{-8pt}{\begingroup\makeatletter\ifx\SetFigFont\undefined
\def\x#1#2#3#4#5#6#7\relax{\def\x{#1#2#3#4#5#6}}%
\expandafter\x\fmtname xxxxxx\relax \def\y{splain}%
\ifx\x\y   
\gdef\SetFigFont#1#2#3{%
  \ifnum #1<17\tiny\else \ifnum #1<20\small\else
  \ifnum #1<24\normalsize\else \ifnum #1<29\large\else
  \ifnum #1<34\Large\else \ifnum #1<41\LARGE\else
     \huge\fi\fi\fi\fi\fi\fi
  \csname #3\endcsname}%
\else
\gdef\SetFigFont#1#2#3{\begingroup
  \count@#1\relax \ifnum 25<\count@\count@25\fi
  \def\x{\endgroup\@setsize\SetFigFont{#2pt}}%
  \expandafter\x
    \csname \romannumeral\the\count@ pt\expandafter\endcsname
    \csname @\romannumeral\the\count@ pt\endcsname
  \csname #3\endcsname}%
\fi
\fi\endgroup
{\renewcommand{\dashlinestretch}{30}
\begin{picture}(5521,1974)(0,-10)
\path(375,1725)(1875,1725)
\path(1755.000,1695.000)(1875.000,1725.000)(1755.000,1755.000)
\path(2355.000,1455.000)(2325.000,1575.000)(2295.000,1455.000)
\path(2325,1575)(2325,375)
\path(525,225)(1875,225)
\path(1755.000,195.000)(1875.000,225.000)(1755.000,255.000)
\path(225,375)(225,1575)
\path(255.000,1455.000)(225.000,1575.000)(195.000,1455.000)
\path(2400,1200)(2402,1197)(2407,1190)
	(2414,1178)(2424,1163)(2435,1145)
	(2446,1126)(2456,1108)(2464,1091)
	(2470,1076)(2474,1063)(2475,1050)
	(2474,1037)(2470,1025)(2465,1012)
	(2457,1000)(2448,987)(2437,975)
	(2427,962)(2418,950)(2410,937)
	(2405,925)(2401,912)(2400,900)
	(2401,887)(2405,874)(2411,859)
	(2419,842)(2429,824)(2440,805)
	(2451,787)(2461,772)(2468,760)
	(2473,753)(2475,750)
\path(300,1200)(302,1197)(307,1190)
	(314,1178)(324,1163)(335,1145)
	(346,1126)(356,1108)(364,1091)
	(370,1076)(374,1063)(375,1050)
	(374,1037)(370,1025)(365,1012)
	(357,1000)(348,987)(337,975)
	(327,962)(318,950)(310,937)
	(305,925)(301,912)(300,900)
	(301,887)(305,874)(311,859)
	(319,842)(329,824)(340,805)
	(351,787)(361,772)(368,760)
	(373,753)(375,750)
\put(0,900){\makebox(0,0)[lb]{\smash{{\mathmode{\ch}}}}}
\put(975,1815){\makebox(0,0)[lb]{\smash{{\mathmode{\Tg}}}}}
\put(1050,315){\makebox(0,0)[lb]{\smash{{\mathmode{\Tg}}}}}
\put(150,150){\makebox(0,0)[lb]{\smash{{\mathmode{\B}}}}}
\put(2025,900){\makebox(0,0)[lb]{\smash{{\mathmode{\b}}}}}
\put(150,1650){\makebox(0,0)[lb]{\smash{{\mathmode{\A}}}}}
\put(2625,1575){\makebox(0,0)[lb]{\smash{{\mathmode{\Udef}}}}}
\put(2625,300){\makebox(0,0)[lb]{\smash{{\mathmode{\Sdef}}}}}
\put(5475,1800){\makebox(0,0)[lb]{\smash{{\mathmode{ }}}}}
\put(5475,0){\makebox(0,0)[lb]{\smash{{\mathmode{ }}}}}
\put(1950,1650){\makebox(0,0)[lb]{\smash{{\mathmode{\U}}}}}
\put(1950,150){\makebox(0,0)[lb]{\smash{{\mathmode{\S}}}}}
\end{picture}
} }
        	\hspace{-1.9mm}
	\end{array}

} \]
in which the left column is the above mentioned formal PBW isomorphism
$\chi$, and the right column is the symmetrization map $\b:S(\frakg)\to
U(\frakg)$, sending an unordered word of length $n$ to the average of the
$n!$ ways of ordering its letters and reading them as a product in
$U(\frakg)$. The map $\b$ is an isomorphism by the honest PBW theorem. The
left-to-right maps ${\mathcal T}_\frakg$ are defined as in
\cite{Bar-Natan:Vassiliev} by contracting copies of the structure
constants tensor, one for each vertex of any given diagram, using the
standard invariant form $(\cdot,\cdot)$ on $\frakg$ (see citations in
section~\ref{TheEdges} below). The maps ${\mathcal T}_\frakg$ seem
to `forget' some information (some high-degree elements on the left
get mapped to 0 on the right no matter what the algebra $\frakg$ is,
see~\cite{Vogel:Structures}), but at least up to degree 12 it is faithful
(for some Lie algebras); see~\cite{Kneissler:Twelve}.

\begin{theorem} \lbl{UnknotTheorem} Conjecture~\ref{UnknotConjecture} is
``true on the level of semi-simple Lie algebras''. Namely,
\[ {\mathcal T}_\frakg\Omega = {\mathcal T}_\frakg\chi^{-1}\ZF(\bigcirc). \]
\end{theorem}

We now formulate our second conjecture.  Let
$\calB'=\Span\left\{\setlength{\unitlength}{0.5\standardunitlength}
	\begin{array}{c}  \hspace{-1.7mm}
        	\raisebox{-2pt}{%
\begingroup\makeatletter\ifx\SetFigFont\undefined
\def\x#1#2#3#4#5#6#7\relax{\def\x{#1#2#3#4#5#6}}%
\expandafter\x\fmtname xxxxxx\relax \def\y{splain}%
\ifx\x\y   
\gdef\SetFigFont#1#2#3{%
  \ifnum #1<17\tiny\else \ifnum #1<20\small\else
  \ifnum #1<24\normalsize\else \ifnum #1<29\large\else
  \ifnum #1<34\Large\else \ifnum #1<41\LARGE\else
     \huge\fi\fi\fi\fi\fi\fi
  \csname #3\endcsname}%
\else
\gdef\SetFigFont#1#2#3{\begingroup
  \count@#1\relax \ifnum 25<\count@\count@25\fi
  \def\x{\endgroup\@setsize\SetFigFont{#2pt}}%
  \expandafter\x
    \csname \romannumeral\the\count@ pt\expandafter\endcsname
    \csname @\romannumeral\the\count@ pt\endcsname
  \csname #3\endcsname}%
\fi
\fi\endgroup
{\renewcommand{\dashlinestretch}{30}
\begin{picture}(920,639)(0,-10)
\put(762,312){\ellipse{300}{300}}
\path(162,462)(462,462)(462,162)
	(162,162)(162,462)
\path(162,162)(12,12)
\path(462,162)(612,12)
\path(462,462)(612,612)
\path(162,462)(12,612)
\path(162,312)(462,312)
\path(762,462)(762,162)
\end{picture}
}
 }
        	\hspace{-1.9mm}
	\end{array}
\right\}/(AS,IHX)$
\lbl{Bdefinition} be the same as $\calB$, only dropping the connectivity
requirement (so that
we also allow
connected components that have no univalent vertices). The space $\calB'$
has two different products, and thus is an algebra in two different ways:
\begin{myitemize}
\item The disjoint union $C_1\udot C_2$ of two Chinese characters $C_{1,2}$
  is again a Chinese character. The obvious bilinear extension of $\udot$
  is a well defined product $\calB'\times\calB'\to\calB'$, which turns
  $\calB'$ into an algebra. For emphasis we will call this algebra $\Bu$.
\item $\calB'$ is isomorphic (as a vector space) to the space
  $\calA'=\Span\left\{\setlength{\unitlength}{0.5\standardunitlength}
	\begin{array}{c}  \hspace{-1.7mm}
        	\raisebox{-2pt}{%
\begingroup\makeatletter\ifx\SetFigFont\undefined
\def\x#1#2#3#4#5#6#7\relax{\def\x{#1#2#3#4#5#6}}%
\expandafter\x\fmtname xxxxxx\relax \def\y{splain}%
\ifx\x\y   
\gdef\SetFigFont#1#2#3{%
  \ifnum #1<17\tiny\else \ifnum #1<20\small\else
  \ifnum #1<24\normalsize\else \ifnum #1<29\large\else
  \ifnum #1<34\Large\else \ifnum #1<41\LARGE\else
     \huge\fi\fi\fi\fi\fi\fi
  \csname #3\endcsname}%
\else
\gdef\SetFigFont#1#2#3{\begingroup
  \count@#1\relax \ifnum 25<\count@\count@25\fi
  \def\x{\endgroup\@setsize\SetFigFont{#2pt}}%
  \expandafter\x
    \csname \romannumeral\the\count@ pt\expandafter\endcsname
    \csname @\romannumeral\the\count@ pt\endcsname
  \csname #3\endcsname}%
\fi
\fi\endgroup
{\renewcommand{\dashlinestretch}{30}
\begin{picture}(1224,409)(0,-10)
\put(837.000,24.500){\arc{450.694}{3.0861}{6.3387}}
\put(612.000,12.000){\arc{300.000}{3.1416}{6.2832}}
\put(237,237){\ellipse{300}{300}}
\path(12,12)(1212,12)
\path(1092.000,-18.000)(1212.000,12.000)(1092.000,42.000)
\path(912,237)(912,12)
\path(87,237)(387,237)
\end{picture}
}
 }
        	\hspace{-1.9mm}
	\end{array}
\right\}/(AS,IHX,STU)$
  of diagrams whose skeleton is a single oriented interval
  (like $\calA$, only that here we also
  allow non-connected diagrams). The isomorphism is the map
  $\chi:\calB'\to\calA'$ that maps a Chinese
  character with $k$ ``legs'' (univalent vertices) to the average
  of the $k!$ ways of arranging them along an oriented interval
  (in~\cite{Bar-Natan:Vassiliev} the sum was used instead of the
  average). $\calA'$ has a well known ``juxtaposition'' product $\times$,
  related to the ``connect sum'' operation on knots:
  \[ \setlength{\unitlength}{0.5\standardunitlength}
	\begin{array}{c}  \hspace{-1.7mm}
        	\raisebox{-2pt}{%
\begingroup\makeatletter\ifx\SetFigFont\undefined
\def\x#1#2#3#4#5#6#7\relax{\def\x{#1#2#3#4#5#6}}%
\expandafter\x\fmtname xxxxxx\relax \def\y{splain}%
\ifx\x\y   
\gdef\SetFigFont#1#2#3{%
  \ifnum #1<17\tiny\else \ifnum #1<20\small\else
  \ifnum #1<24\normalsize\else \ifnum #1<29\large\else
  \ifnum #1<34\Large\else \ifnum #1<41\LARGE\else
     \huge\fi\fi\fi\fi\fi\fi
  \csname #3\endcsname}%
\else
\gdef\SetFigFont#1#2#3{\begingroup
  \count@#1\relax \ifnum 25<\count@\count@25\fi
  \def\x{\endgroup\@setsize\SetFigFont{#2pt}}%
  \expandafter\x
    \csname \romannumeral\the\count@ pt\expandafter\endcsname
    \csname @\romannumeral\the\count@ pt\endcsname
  \csname #3\endcsname}%
\fi
\fi\endgroup
{\renewcommand{\dashlinestretch}{30}
\begin{picture}(774,186)(0,-10)
\put(312.000,12.000){\arc{300.000}{3.1416}{6.2832}}
\put(462.000,12.000){\arc{300.000}{3.1416}{6.2832}}
\path(12,12)(762,12)
\path(642.000,-18.000)(762.000,12.000)(642.000,42.000)
\end{picture}
}
 }
        	\hspace{-1.9mm}
	\end{array}
\times\ \setlength{\unitlength}{0.5\standardunitlength}
	\begin{array}{c}  \hspace{-1.7mm}
        	\raisebox{-2pt}{%
\begingroup\makeatletter\ifx\SetFigFont\undefined
\def\x#1#2#3#4#5#6#7\relax{\def\x{#1#2#3#4#5#6}}%
\expandafter\x\fmtname xxxxxx\relax \def\y{splain}%
\ifx\x\y   
\gdef\SetFigFont#1#2#3{%
  \ifnum #1<17\tiny\else \ifnum #1<20\small\else
  \ifnum #1<24\normalsize\else \ifnum #1<29\large\else
  \ifnum #1<34\Large\else \ifnum #1<41\LARGE\else
     \huge\fi\fi\fi\fi\fi\fi
  \csname #3\endcsname}%
\else
\gdef\SetFigFont#1#2#3{\begingroup
  \count@#1\relax \ifnum 25<\count@\count@25\fi
  \def\x{\endgroup\@setsize\SetFigFont{#2pt}}%
  \expandafter\x
    \csname \romannumeral\the\count@ pt\expandafter\endcsname
    \csname @\romannumeral\the\count@ pt\endcsname
  \csname #3\endcsname}%
\fi
\fi\endgroup
{\renewcommand{\dashlinestretch}{30}
\begin{picture}(1224,409)(0,-10)
\put(837.000,24.500){\arc{450.694}{3.0861}{6.3387}}
\put(612.000,12.000){\arc{300.000}{3.1416}{6.2832}}
\put(237,237){\ellipse{300}{300}}
\path(12,12)(1212,12)
\path(1092.000,-18.000)(1212.000,12.000)(1092.000,42.000)
\path(912,237)(912,12)
\path(87,237)(387,237)
\end{picture}
}
 }
        	\hspace{-1.9mm}
	\end{array}
=
    \setlength{\unitlength}{0.5\standardunitlength}
	\begin{array}{c}  \hspace{-1.7mm}
        	\raisebox{-2pt}{%
\begingroup\makeatletter\ifx\SetFigFont\undefined
\def\x#1#2#3#4#5#6#7\relax{\def\x{#1#2#3#4#5#6}}%
\expandafter\x\fmtname xxxxxx\relax \def\y{splain}%
\ifx\x\y   
\gdef\SetFigFont#1#2#3{%
  \ifnum #1<17\tiny\else \ifnum #1<20\small\else
  \ifnum #1<24\normalsize\else \ifnum #1<29\large\else
  \ifnum #1<34\Large\else \ifnum #1<41\LARGE\else
     \huge\fi\fi\fi\fi\fi\fi
  \csname #3\endcsname}%
\else
\gdef\SetFigFont#1#2#3{\begingroup
  \count@#1\relax \ifnum 25<\count@\count@25\fi
  \def\x{\endgroup\@setsize\SetFigFont{#2pt}}%
  \expandafter\x
    \csname \romannumeral\the\count@ pt\expandafter\endcsname
    \csname @\romannumeral\the\count@ pt\endcsname
  \csname #3\endcsname}%
\fi
\fi\endgroup
{\renewcommand{\dashlinestretch}{30}
\begin{picture}(1825,409)(0,-10)
\put(1438.000,24.500){\arc{450.694}{3.0861}{6.3387}}
\put(1213.000,12.000){\arc{300.000}{3.1416}{6.2832}}
\put(312.000,12.000){\arc{300.000}{3.1416}{6.2832}}
\put(462.000,12.000){\arc{300.000}{3.1416}{6.2832}}
\put(838,237){\ellipse{300}{300}}
\path(12,12)(1813,12)
\path(1693.000,-18.000)(1813.000,12.000)(1693.000,42.000)
\path(1513,237)(1513,12)
\path(688,237)(988,237)
\end{picture}
}
 }
        	\hspace{-1.9mm}
	\end{array}
.
  \]
  The algebra structure on $\calA'$ defines another algebra structure on
  $\calB'$. For emphasis we will call this algebra $\Bt$.
\end{myitemize}

As before, $\calA'$ is graded by half the number of trivalent
vertices in a diagram, $\calB'$ is graded by half the total number of
vertices in a diagram, and the isomorphism $\chi$ as well as the two
products respect these gradings.

\begin{definition} \lbl{ChatDef} If $C$ is a Chinese character, let
$\hat{C}:\calB'\to\calB'$ be the operator defined by
\[ \hat{C}(C')=\begin{cases}
    0 & \text{if $C$ has more legs than $C'$,} \\
    \parbox{2.6in}{
      the sum of all ways of gluing all the legs of $C$ to some (or all)
      legs of $C'$
    }\quad & \text{otherwise.}
  \end{cases}
\]
For example,
\[ \widehat{\omega_4}(\omega_2)=0; \qquad
   \widehat{\omega_2}(\omega_4)=
     8\printname{SideGluing}
	\setlength{\unitlength}{0.5\standardunitlength}
	\begin{array}{c}  \hspace{-1.7mm}
        	\raisebox{-8pt}{%
\begingroup\makeatletter\ifx\SetFigFont\undefined
\def\x#1#2#3#4#5#6#7\relax{\def\x{#1#2#3#4#5#6}}%
\expandafter\x\fmtname xxxxxx\relax \def\y{splain}%
\ifx\x\y   
\gdef\SetFigFont#1#2#3{%
  \ifnum #1<17\tiny\else \ifnum #1<20\small\else
  \ifnum #1<24\normalsize\else \ifnum #1<29\large\else
  \ifnum #1<34\Large\else \ifnum #1<41\LARGE\else
     \huge\fi\fi\fi\fi\fi\fi
  \csname #3\endcsname}%
\else
\gdef\SetFigFont#1#2#3{\begingroup
  \count@#1\relax \ifnum 25<\count@\count@25\fi
  \def\x{\endgroup\@setsize\SetFigFont{#2pt}}%
  \expandafter\x
    \csname \romannumeral\the\count@ pt\expandafter\endcsname
    \csname @\romannumeral\the\count@ pt\endcsname
  \csname #3\endcsname}%
\fi
\fi\endgroup
{\renewcommand{\dashlinestretch}{30}
\begin{picture}(695,639)(0,-10)
\put(83,312){\ellipse{150}{150}}
\path(533,462)(233,462)(233,162)
	(533,162)(533,462)
\path(683,612)(533,462)
\path(533,162)(683,12)
\path(233,462)	(197.611,520.594)
	(158.000,537.000)

\path(158,537)	(109.137,509.411)
	(83.000,462.000)

\path(83,462)	(78.642,435.086)
	(83.000,387.000)

\path(83,237)	(78.650,188.914)
	(83.000,162.000)

\path(83,162)	(109.137,114.589)
	(158.000,87.000)

\path(158,87)	(197.611,103.406)
	(233.000,162.000)

\end{picture}
}
 }
        	\hspace{-1.9mm}
	\end{array}
+4\printname{DiagonalGluing}
	\setlength{\unitlength}{0.5\standardunitlength}
	\begin{array}{c}  \hspace{-1.7mm}
        	\raisebox{-8pt}{%
\begingroup\makeatletter\ifx\SetFigFont\undefined
\def\x#1#2#3#4#5#6#7\relax{\def\x{#1#2#3#4#5#6}}%
\expandafter\x\fmtname xxxxxx\relax \def\y{splain}%
\ifx\x\y   
\gdef\SetFigFont#1#2#3{%
  \ifnum #1<17\tiny\else \ifnum #1<20\small\else
  \ifnum #1<24\normalsize\else \ifnum #1<29\large\else
  \ifnum #1<34\Large\else \ifnum #1<41\LARGE\else
     \huge\fi\fi\fi\fi\fi\fi
  \csname #3\endcsname}%
\else
\gdef\SetFigFont#1#2#3{\begingroup
  \count@#1\relax \ifnum 25<\count@\count@25\fi
  \def\x{\endgroup\@setsize\SetFigFont{#2pt}}%
  \expandafter\x
    \csname \romannumeral\the\count@ pt\expandafter\endcsname
    \csname @\romannumeral\the\count@ pt\endcsname
  \csname #3\endcsname}%
\fi
\fi\endgroup
{\renewcommand{\dashlinestretch}{30}
\begin{picture}(695,658)(0,-10)
\put(83,312){\ellipse{150}{150}}
\path(533,462)(233,462)(233,162)
	(533,162)(533,462)
\path(533,162)(683,12)
\path(233,462)(83,612)
\path(533,462)	(591.594,497.389)
	(608.000,537.000)

\path(608,537)	(580.411,585.862)
	(533.000,612.000)

\path(533,612)	(465.104,629.415)
	(424.940,633.769)
	(383.000,635.220)
	(341.060,633.769)
	(300.896,629.415)
	(233.000,612.000)

\path(233,612)	(190.345,586.760)
	(146.390,548.610)
	(108.240,504.655)
	(83.000,462.000)

\path(83,462)	(78.643,435.086)
	(83.000,387.000)

\path(83,237)	(78.650,188.914)
	(83.000,162.000)

\path(83,162)	(109.138,114.589)
	(158.000,87.000)

\path(158,87)	(197.611,103.406)
	(233.000,162.000)

\end{picture}
}
 }
        	\hspace{-1.9mm}
	\end{array}
.
\]
If $C$ has $k$ legs and total degree $m$, then $\hat{C}$ is an operator
of degree $m-k$. By linear extension, we find that every $C\in\calB'$
defines an operator $\hat{C}:\calB'\to\calB'$, and in fact, even infinite
linear combinations of Chinese characters with an {\em increasing} number
of legs define operators $\calB'\to\calB'$.
\end{definition}

As $\Omega$ is made of wheels, we call the action of the (degree
$0$) operator $\hat{\Omega}$ \lbl{Ohdefinition} ``wheeling''. As
$\Omega$ begins with $1$, the wheeling map is invertible. We argue
below that $\hat{\Omega}$ is a diagrammatic analog of the Duflo
isomorphism ${S}^\frakg(\frakg)\to{S}^\frakg(\frakg)$
(see~\cite{Duflo:Resolubles} and see below). The Duflo isomorphism
intertwines the two algebra structures that ${S}^\frakg(\frakg)$ has:
the structure it inherits from the symmetric algebra and the structure
it inherits from ${U}^\frakg(\frakg)$ via the PBW isomorphism.
One may hope that $\hat{\Omega}$ has the parallel property:

\begin{conjecture} \lbl{WheelingConjecture} (Wheeling\footnote{Conjectured
independently by Deligne~\cite{Deligne:Letter}.}) Wheeling intertwines
the two products on Chinese characters. More precisely, the map
$\hat{\Omega}:\Bu\to\Bt$ is an algebra isomorphism.
\end{conjecture}

There are several good reasons to hope that
Conjecture~\ref{WheelingConjecture} is true. If it is true, one would
be able to use it along with Conjecture~\ref{UnknotConjecture} and
known properties of the Kontsevich integral (such as its behavior the
operations of change of framing, connected sum, and taking the parallel of
a component as in~\cite{LeMurakami:Parallel}) to get explicit formulas
for the Kontsevich integral of several other knots and links. Note
that change of framing and connect sum act on the Kontsevich integral
multiplicatively using the product in $\calA$, but the conjectured formula
we have for the Kontsevich integral of the unknot is in $\calB$. Using
Conjecture~\ref{WheelingConjecture} it should be possible to perform
all operations in $\calB$.

Perhaps a more important reason is that in essence, $\calA$ and $\calB$
capture that part of the information about $U(\frakg)$ and $S(\frakg)$
that can be described entirely in terms of the bracket and the structure
constants. Thus a proof of Conjecture~\ref{WheelingConjecture} would yield an
elementary proof of the intertwining property of the Duflo isomorphism,
whose current proofs use representation theory and are quite involved. We
feel that the knowledge missing to give an elementary proof of the
intertwining property of the Duflo isomorphism is the same knowledge
that is missing for giving a proof of the Kashiwara-Vergne conjecture
(\cite{KashiwaraVergne:CampbellHausdorff}).

\begin{theorem} \lbl{WheelingTheorem} Conjecture~\ref{WheelingConjecture}
is ``true on the level of semi-simple Lie algebras''. A precise statement is
in Proposition~\ref{PreciseWheelingTheorem} and the remark following it.
\end{theorem}

\begin{remark} As semi-simple Lie algebras ``see'' all of the
Vassiliev theory at least up to degree~12 \cite{Bar-Natan:Vassiliev,
Kneissler:Twelve}, Theorems~\ref{UnknotTheorem} and ~\ref{WheelingTheorem}
imply Conjectures~\ref{UnknotConjecture} and~\ref{WheelingConjecture} up to that
degree. It should be noted that semi-simple Lie algebras do not ``see''
the whole Vassiliev theory at high degrees, see ~\cite{Vogel:Structures}.
\end{remark}

\begin{remark} As the Duflo isomorphism has no known elementary proof, the
Lie algebra techniques used in this paper are unlikely to give full proofs
of Conjectures~\ref{UnknotConjecture} and~\ref{WheelingConjecture}.
\end{remark}

\begin{remark} We've chosen to work over the complex numbers to allow for
some analytical arguments below. The rationality of the Kontsevich
integral~\cite{LeMurakami:Universal} and the uniform classification of
semi-simple Lie algebras over fields of characteristic 0 implies that
Conjectures~\ref{UnknotConjecture} and~\ref{WheelingConjecture} and
Theorems~\ref{UnknotTheorem} and~\ref{WheelingTheorem} are independent
of the (characteristic 0) ground field.
\end{remark}

\subsection{The plan} Theorem~\ref{UnknotTheorem} and
Theorem~\ref{WheelingTheorem} both follow from a delicate assembly of
widely known facts about Lie algebras and related objects; the main
novelty in this paper is the realization that these known facts can
be brought together and used to prove Theorems~\ref{UnknotTheorem}
and~\ref{WheelingTheorem} and make Conjectures~\ref{UnknotConjecture}
and~\ref{WheelingConjecture}. The facts we use about Lie-algebras amount to the
commutativity of a certain monstrous diagram.  In Section~\ref{diagram}
below we will explain everything that appears in that diagram,
prove its commutativity, and prove Theorem~\ref{WheelingTheorem}. In
Section~\ref{proof} we will show how that commutativity implies
Theorem~\ref{UnknotTheorem} as well. We conclude this introductory
section with a picture of the monster itself:

\begin{theorem} \lbl{CommutativityTheorem} (definitions and proof in
Section~\ref{diagram})
The following monster diagram is commutative:
\[ \printname{CommutativeDiagram}
	\setlength{\unitlength}{0.8\standardunitlength}
	\begin{array}{c}  \hspace{-1.7mm}
        	\raisebox{-8pt}{\begingroup\makeatletter\ifx\SetFigFont\undefined
\def\x#1#2#3#4#5#6#7\relax{\def\x{#1#2#3#4#5#6}}%
\expandafter\x\fmtname xxxxxx\relax \def\y{splain}%
\ifx\x\y   
\gdef\SetFigFont#1#2#3{%
  \ifnum #1<17\tiny\else \ifnum #1<20\small\else
  \ifnum #1<24\normalsize\else \ifnum #1<29\large\else
  \ifnum #1<34\Large\else \ifnum #1<41\LARGE\else
     \huge\fi\fi\fi\fi\fi\fi
  \csname #3\endcsname}%
\else
\gdef\SetFigFont#1#2#3{\begingroup
  \count@#1\relax \ifnum 25<\count@\count@25\fi
  \def\x{\endgroup\@setsize\SetFigFont{#2pt}}%
  \expandafter\x
    \csname \romannumeral\the\count@ pt\expandafter\endcsname
    \csname @\romannumeral\the\count@ pt\endcsname
  \csname #3\endcsname}%
\fi
\fi\endgroup
{\renewcommand{\dashlinestretch}{30}
\begin{picture}(5787,4659)(0,-10)
\path(375,4461)(375,3261)
\path(345.000,3381.000)(375.000,3261.000)(405.000,3381.000)
\path(525,3111)(2025,3111)
\path(1905.000,3081.000)(2025.000,3111.000)(1905.000,3141.000)
\path(675,4461)(2175,3261)
\path(2062.555,3312.537)(2175.000,3261.000)(2100.037,3359.389)
\path(2475,2961)(2475,1761)
\path(2505.000,2841.000)(2475.000,2961.000)(2445.000,2841.000)
\path(525,1461)(1275,261)
\path(614.040,1375.140)(525.000,1461.000)(563.160,1343.340)
\path(1575,111)(3075,111)
\path(2955.000,81.000)(3075.000,111.000)(2955.000,141.000)
\path(675,1611)(2025,1611)
\path(1905.000,1581.000)(2025.000,1611.000)(1905.000,1641.000)
\path(2625,1461)(3525,261)
\path(2721.000,1383.000)(2625.000,1461.000)(2673.000,1347.000)
\path(4725,1461)(5625,261)
\path(4821.000,1383.000)(4725.000,1461.000)(4773.000,1347.000)
\path(4950,2961)(5775,261)
\path(5013.757,2855.004)(4950.000,2961.000)(4956.376,2837.471)
\path(3075,3111)(4275,3111)
\path(4155.000,3081.000)(4275.000,3111.000)(4155.000,3141.000)
\path(3075,1636)(4125,1636)
\path(3075,1586)(4125,1586)
\path(4125,135)(5175,135)
\path(4125,85)(5175,85)
\path(375,1761)(375,2961)
\path(405.000,2841.000)(375.000,2961.000)(345.000,2841.000)
\path(2550,2586)	(2585.449,2535.329)
	(2608.594,2495.783)
	(2625.000,2436.000)

\path(2625,2436)	(2613.281,2399.079)
	(2587.500,2361.000)
	(2561.719,2322.921)
	(2550.000,2286.000)

\path(2550,2286)	(2566.406,2226.217)
	(2589.551,2186.671)
	(2625.000,2136.000)

\path(5100,1161)	(5161.662,1112.399)
	(5204.872,1073.543)
	(5250.000,1011.000)

\path(5250,1011)	(5254.354,977.052)
	(5250.000,936.000)
	(5245.646,894.948)
	(5250.000,861.000)

\path(5250,861)	(5295.128,798.457)
	(5338.338,759.601)
	(5400.000,711.000)

\path(3375,3186)	(3425.671,3221.449)
	(3465.217,3244.594)
	(3525.000,3261.000)

\path(3525,3261)	(3561.921,3249.281)
	(3600.000,3223.500)
	(3638.079,3197.719)
	(3675.000,3186.000)

\path(3675,3186)	(3734.783,3202.406)
	(3774.329,3225.551)
	(3825.000,3261.000)

\path(3000,1161)	(3061.662,1112.399)
	(3104.872,1073.543)
	(3150.000,1011.000)

\path(3150,1011)	(3154.354,977.052)
	(3150.000,936.000)
	(3145.646,894.948)
	(3150.000,861.000)

\path(3150,861)	(3195.128,798.457)
	(3238.338,759.601)
	(3300.000,711.000)

\path(825,1161)	(886.662,1112.399)
	(929.872,1073.543)
	(975.000,1011.000)

\path(975,1011)	(979.354,977.052)
	(975.000,936.000)
	(970.646,894.948)
	(975.000,861.000)

\path(975,861)	(1020.128,798.457)
	(1063.338,759.601)
	(1125.000,711.000)

\path(5325,1986)	(5386.662,1937.399)
	(5429.872,1898.543)
	(5475.000,1836.000)

\path(5475,1836)	(5482.491,1784.807)
	(5477.561,1722.780)
	(5471.351,1661.113)
	(5475.000,1611.000)

\path(5475,1611)	(5490.931,1572.124)
	(5518.898,1525.181)
	(5562.415,1464.898)
	(5591.105,1428.105)
	(5625.000,1386.000)

\path(450,2586)	(485.449,2535.329)
	(508.594,2495.783)
	(525.000,2436.000)

\path(525,2436)	(513.281,2399.079)
	(487.500,2361.000)
	(461.719,2322.921)
	(450.000,2286.000)

\path(450,2286)	(466.406,2226.217)
	(489.551,2186.671)
	(525.000,2136.000)

\put(375,4536){\makebox(0,0)[lb]{\smash{{\mathmode{\KF}}}}}
\put(150,2286){\makebox(0,0)[lb]{\smash{{\mathmode{\ch}}}}}
\put(300,3036){\makebox(0,0)[lb]{\smash{{\mathmode{\A}}}}}
\put(1500,3861){\makebox(0,0)[lb]{\smash{{\mathmode{\RTg}}}}}
\put(1125,3201){\makebox(0,0)[lb]{\smash{{\mathmode{\Tg}}}}}
\put(1200,1701){\makebox(0,0)[lb]{\smash{{\mathmode{\Tg}}}}}
\put(2175,201){\makebox(0,0)[lb]{\smash{{\mathmode{\Tg}}}}}
\put(0,3786){\makebox(0,0)[lb]{\smash{{\mathmode{\ZF}}}}}
\put(675,711){\makebox(0,0)[lb]{\smash{{\mathmode{\Oh}}}}}
\put(4875,711){\makebox(0,0)[lb]{\smash{{\mathmode{\Sg}}}}}
\put(525,4161){\makebox(0,0)[lb]{\smash{{\mathmode{\Ca}}}}}
\put(600,2736){\makebox(0,0)[lb]{\smash{{\mathmode{\Cb}}}}}
\put(2700,2661){\makebox(0,0)[lb]{\smash{{\mathmode{\Cc}}}}}
\put(825,1236){\makebox(0,0)[lb]{\smash{{\mathmode{\Cd}}}}}
\put(3000,1236){\makebox(0,0)[lb]{\smash{{\mathmode{\Ce}}}}}
\put(300,1536){\makebox(0,0)[lb]{\smash{{\mathmode{\Bt}}}}}
\put(1275,36){\makebox(0,0)[lb]{\smash{{\mathmode{\Bu}}}}}
\put(2025,3036){\makebox(0,0)[lb]{\smash{{\mathmode{\Ug}}}}}
\put(2025,1536){\makebox(0,0)[lb]{\smash{{\mathmode{\St}}}}}
\put(3075,36){\makebox(0,0)[lb]{\smash{{\mathmode{\Su}}}}}
\put(4200,1536){\makebox(0,0)[lb]{\smash{{\mathmode{\Pg}}}}}
\put(5250,36){\makebox(0,0)[lb]{\smash{{\mathmode{\Pg}}}}}
\put(4350,3036){\makebox(0,0)[lb]{\smash{{\mathmode{\Ph}}}}}
\put(3450,2886){\makebox(0,0)[lb]{\smash{{\mathmode{\psig}}}}}
\put(2175,2286){\makebox(0,0)[lb]{\smash{{\mathmode{\b}}}}}
\put(5550,1761){\makebox(0,0)[lb]{\smash{{\mathmode{\ig}}}}}
\put(2250,711){\makebox(0,0)[lb]{\smash{{\mathmode{\Dj}}}}}
\end{picture}
} }
        	\hspace{-1.9mm}
	\end{array}
 \]
\end{theorem}

\begin{remark} Our two conjectures ought to be related---one talks
about $\Omega$, and another is about an operator $\hat{\Omega}$
made out of $\Omega$, and the proofs of Theorems~\ref{UnknotTheorem}
and~\ref{WheelingTheorem} both use the Duflo map ($\Dj$ in the above
diagram). But looking more closely at the proofs below, the relationship
seems to disappear. The proof of Theorem~\ref{WheelingTheorem} uses
only the commutativity of the face labeled $\Cd$, while the proof of
Theorem~\ref{UnknotTheorem} uses the commutativity of all faces but
$\Cd$. No further relations between the conjectures are seen in the
proofs of our theorems. We are still missing the deep relation that
ought to exist between `Wheels' and `Wheeling'. Why is it that the same
strange combination of Chinese characters $\Omega$ plays a role in these
two seemingly unrelated affairs?
\end{remark}

\subsection{Postscript} According to
Kontsevich~\cite{Kontsevich:DeformationQuantization},
Conjecture~\ref{WheelingConjecture} seems to follow from the results he proves
in Section~8.3 of that paper, but a full proof of the conjecture is not
given there. \cite{LeThurston:Unknot} have shown that
Conjecture~\ref{WheelingConjecture} implies
Conjecture~\ref{UnknotConjecture}, but
unfortunately their proof does not shed light on the fundamental
relationship that ought to exist between the two conjectures.

\subsection{Acknowledgement} Much of this work was done when the
four of us were visiting \Arhus, Denmark, for a special semester on
geometry and physics, in August 1995. We wish to thank the organizers,
J.~Dupont, H.~Pedersen, A.~Swann and especially J.~Andersen for their
hospitality and for the stimulating atmosphere they created. We wish
to thank the Institute for Advanced Studies for their hospitality,
and P.~Deligne for listening to our thoughts and sharing his. His
letter~\cite{Deligne:Letter} introduced us to the Duflo isomorphism;
initially our proofs relied more heavily on the Kirillov character
formula. A.~Others made some very valuable suggestions; we thank them
and also thank J.~Birman, A.~Haviv, A.~Joseph, G.~Perets, J.~D.~Rogawski,
M.~Vergne and S.~Willerton for additional remarks and suggestions.

\section{The monster diagram} \lbl{diagram}

\subsection{The vertices} Let $\frakg_{\mathbf R}$ be the (semi-simple) 
Lie-algebra of some compact Lie group $G$, let $\frakg=\frakg_{\mathbf R}
\otimes{\mathbf C}$, let $\frakh\subset i\frakg_{\mathbf R}$ be a Cartan
subalgebra of $\frakg$, and let $W$ be the Weyl group of $\frakh$ in
$\frakg$. Let $\Delta_+\subset\frakh^\star$ be a set of positive
roots of $\frakg$, and let $\rho\in i\frakg_{\mathbf R}^\star$ be half
the sum of the positive roots.  Let $\hbar$ be an indeterminate, and
let ${\mathbf C}[[\hbar]]$ be the ring of formal power series in $\hbar$
with coefficients in ${\mathbf C}$.

\begin{myitemize}

\item $\KF$ is the set of all framed knots in ${\mathbf R}^3$.

\item $\A$ is the algebra of not-necessarily-connected chord diagrams,
  as in page~\pageref{Adefinition}.

\item $\Bt$ and $\Bu$ denote the space of Chinese characters (allowing
connected components that have no univalent vertices), as in
page~\pageref{Bdefinition}, taken with its two algebra structures.

\item $\Ug$ is the $\frakg$-invariant part of the universal enveloping
algebra $U(\frakg)$ of $\frakg$, with the coefficient ring extended to be
${\mathbf C}[[\hbar]]$.

\item $\St$ \lbl{Stdefinition} and $\Su$ denote the $\frakg$-invariant
part of the symmetric algebra $S(\frakg)$ of $\frakg$, with the
coefficient ring extended to be ${\mathbf C}[[\hbar]]$. In $\Su$ we take
the algebra structure induced from the natural algebra structure of the
symmetric algebra. In $\St$ we take the algebra structure induced from
the algebra structure of $\Ug$ by the symmetrization map $\b:\St\to\Ug$,
which is a linear isomorphism by the Poincare-Birkhoff-Witt theorem.

\item $\Ph$ is the space of Weyl-invariant polynomial functions on
$\frakh^\star$, with coefficients in ${\mathbf C}[[\hbar]]$.

\item $\Pg$ is the space of ad-invariant polynomial functions on
$\frakg^\star$, with coefficients in ${\mathbf C}[[\hbar]]$.

\end{myitemize}

\subsection{The edges} \lbl{TheEdges}

\begin{myitemize}

\item $\ZF$ is the framed version of the Kontsevich integral for knots as
defined in~\cite{LeMurakami:Universal}. A simpler (and equal) definition for
a framed knot $K$ is
\[
  \ZF(K)=e^{\Theta\cdot\text{writhe}(K)/2}
  \cdot S\left(\tilde{Z}(K)\right)
  \in\calA\subset\A,
\]
where $\Theta$ is the chord diagram $\setlength{\unitlength}{0.5\standardunitlength}
	\begin{array}{c}  \hspace{-1.7mm}
        	\raisebox{-2pt}{%
\begingroup\makeatletter\ifx\SetFigFont\undefined
\def\x#1#2#3#4#5#6#7\relax{\def\x{#1#2#3#4#5#6}}%
\expandafter\x\fmtname xxxxxx\relax \def\y{splain}%
\ifx\x\y   
\gdef\SetFigFont#1#2#3{%
  \ifnum #1<17\tiny\else \ifnum #1<20\small\else
  \ifnum #1<24\normalsize\else \ifnum #1<29\large\else
  \ifnum #1<34\Large\else \ifnum #1<41\LARGE\else
     \huge\fi\fi\fi\fi\fi\fi
  \csname #3\endcsname}%
\else
\gdef\SetFigFont#1#2#3{\begingroup
  \count@#1\relax \ifnum 25<\count@\count@25\fi
  \def\x{\endgroup\@setsize\SetFigFont{#2pt}}%
  \expandafter\x
    \csname \romannumeral\the\count@ pt\expandafter\endcsname
    \csname @\romannumeral\the\count@ pt\endcsname
  \csname #3\endcsname}%
\fi
\fi\endgroup
{\renewcommand{\dashlinestretch}{30}
\begin{picture}(1074,336)(0,-10)
\put(537.000,12.000){\arc{600.000}{3.1416}{6.2832}}
\path(12,12)(1062,12)
\path(942.000,-18.000)(1062.000,12.000)(942.000,42.000)
\end{picture}
}
 }
        	\hspace{-1.9mm}
	\end{array}
$, $S$ is the
standard algebra map $\calA^r=\calA/<\Theta>\to\calA$ defined by mapping
$\Theta$ to $0$ and leaving all other primitives of $\calA$ in place, and
$\tilde{Z}$ is the Kontsevich integral as in~\cite{Kontsevich:Vassiliev}.

\item $\chi$ is the symmetrization map $\Bt\to\A$, as on
page~\pageref{chidefinition}. It is an algebra isomorphism
by~\cite{Bar-Natan:Vassiliev} and the definition of $\times$.

\item $\Oh$ is the wheeling map as in page~\pageref{Ohdefinition}. We
argue that it should be an algebra (iso-)morphism
(Conjecture~\ref{WheelingConjecture}).

\item $\RTg$ denotes the Reshetikhin-Turaev knot invariant
associated with the Lie algebra $\frakg$~\cite{Reshetikhin:QUEA,
Reshetikhin:Quasitriangle, ReshetikhinTuraev:Ribbon, Turaev:Invariants}.

\item $\Tg$ (in all three instances) is the usual
``diagrams to Lie algebras'' map, as in~\cite[Section~2.4 and
exercise~5.1]{Bar-Natan:Vassiliev}. The only variation we make
is that we multiply the image of a degree $m$ element of $\A$
(or $\Bt$ or $\Bu$) by $\hbar^m$. In the construction of $\Tg$
an invariant bilinear form on $\frakg$ is needed.  We use
the standard form $(\cdot,\cdot)$ used in~\cite{ReshetikhinTuraev:Ribbon}
and in~\cite[Appendix]{ChariPressley:QuantumGroups}. See
also~\cite[Chapter~2]{Kac:InfiniteDimensionalLieAlgebras}.

\item The isomorphism $\b$ was already discussed when $\St$ was defined on
page~\pageref{Stdefinition}.

\item The definition of the ``Duflo map'' $\Dj$ requires some
preliminaries. If $V$ is a vector space, there is an algebra map
$D:P(V)\to\Diff(V^\star)$ between the algebra $P(V)$ of polynomial
functions on $V$ and the algebra $\Diff(V^\star)$ of constant coefficients
differential operators on the symmetric algebra $S(V)$. $D$ is defined on
generators as follows: If $\alpha\in V^\star$ is a degree 1 polynomial
on $V$, set $D(\alpha)(v)=\alpha(v)$ for $v\in V\subset S(V)$, and
extend $D(\alpha)$ to all of $S(V)$ using the Leibnitz law. A different
(but less precise) way of defining $D$ is via the Fourier transform:
Identify $S(V)$ with the space of functions on $V^\star$. A polynomial
function on $V$ becomes a differential operator on $V^\star$ after
taking the Fourier transform, and this defines our map $D$. Either way,
if $j\in P(V)$ is homogeneous of degree $k$, the differential operator
$D(j)$ lowers degrees by $k$ and thus vanishes on the low degrees of
$S(V)$. Hence $D(j)$ makes sense even when $j$ is a power series instead
of a polynomial. This definition has a natural extension to the case
when the spaces involved are extended by ${\mathbf C}[[\hbar]]$, or even
${\mathbf C}(\!(\hbar)\!)$, the algebra of Laurent polynomials in $\hbar$.

Now use this definition of $D$ with $V=\frakg$ to
define the Duflo map $\Dj$, where $j_\frakg(X)$ is defined for $X\in\frakg$ by
\[ j_\frakg(X) =
  \det\left(\frac{\sinh\ad X/2}{\ad X/2}\right).
\]
The square root $j_\frakg^{1/2}$ of
$j_\frakg$ is defined as in~\cite{Duflo:Resolubles}
or~\cite[Section~8.2]{BerlineGetzlerVergne:DiracOperators}, and is a power
series in $X$ that begins with $1$. We note that by Kirillov's formula
for the character of the trivial representation (see e.g.~\cite[Theorem
8.4 with $\lambda=i\rho$]{BerlineGetzlerVergne:DiracOperators}),
$j_\frakg^{1/2}$ is the Fourier transform of the symplectic
measure on $M_{i\rho}$, where $M_{i\rho}$ is the co-adjoint
orbit of $i\rho$ in $\frakg_{\mathbf R}^\star$ (see
e.g.~\cite[Section~7.5]{BerlineGetzlerVergne:DiracOperators}):
\begin{equation} \lbl{FourierTransform}
  j_\frakg^{1/2}(X) = \int_{r\in M_{i\rho}} e^{ir(X)}dr.
\end{equation}
(We consider the symplectic measure as a measure on
$\frakg_{\mathbf R}^\star$, whose support is the subset $M_{i\rho}$
of $\frakg_{\mathbf R}^\star$. Its Fourier transform is a function on
$\frakg_{\mathbf R}$ that can be computed via integration on the support
$M_{i\rho}\subset\frakg_{\mathbf R}^\star$ of the symplectic measure.)
Duflo~\cite[th\'eor\`eme~V.2]{Duflo:Resolubles} proved that $\Dj$ is an
algebra isomorphism.

\item $\psig$ is the Harish-Chandra isomorphism $U(\frakg)^\frakg\to
P(\frakh^\star)^W$ extended by $\hbar$. Using the representation theory
of $\frakg$, it is defined as follows.  If $z$ is in $U(\frakg)^\frakg$
and $\lambda\in\frakh^\star$ is a positive integral weight, we set
$\psig(z)(\lambda)$ to be the scalar by which $z$ acts on the irreducible
representation of $\frakg$ whose heighest weight is $\lambda-\rho$. It
is well known (see e.g.~\cite[Section~23.3]{Humphreys:LieAlgebras})
that this partial definition of $\psig(z)$ extends uniquely to a
Weyl-invariant polynomial (also denoted $\psig(z)$) on $\frakh^\star$,
and that the resulting map $\psig:U(\frakg)^\frakg\to P(\frakh^\star)^W$
is an isomorphism.

\item The two equalities at the lower right quarter of the monster diagram
need no explanation. We note though that if the space of polynomials $\Pg$
is endowed with its obvious algebra structure, only the lower equality is
in fact an equality of algebras.

\item $\ig$ is the restriction map induced by the identification of
$\frakh^\star$ with a subspace of $\frakg^\star$ defined using the
form $(\cdot,\cdot)$ of $\frakg$. The map $\ig$ is an isomorphism by
Chevalley's theorem (see e.g.~\cite[Section~23.1]{Humphreys:LieAlgebras}
and~\cite[Section~VI-2]{BrockerTomDieck:Representations}).

\item $\Sg$ is the extension by $\hbar$ of an integral operator. If
$p(\lambda)$ is an invariant polynomial of $\lambda\in\frakg^\star$, then
\[ \Sg(p)(\lambda)=\int_{r\in M_{i\rho}}p(\lambda-ir)dr. \]
$\Sg$ can also be viewed as a convolution operator (with a measure
concentrated on $M_\rho$), and like all convolution operators, it maps
polynomials to polynomials.

\end{myitemize}

\subsection{The faces}

\begin{myitemize}

\item The commutativity of the face labeled $\Ca$
was proven by Kassel~\cite{Kassel:QuantumGroups} and
Le and Murakami~\cite{LeMurakami:Universal} following
Drinfel'd~\cite{Drinfeld:QuasiHopf, Drinfeld:GalQQ}. We comment that
it is this commutativity that makes the notion of ``canonical Vassiliev
invariants''~\cite{Bar-NatanGaroufalidis:MMR} interesting.

\item The commutativity of the face labeled $\Cb$ is immediate from the
definitions, and was already noted in~\cite{Bar-Natan:Vassiliev}.

\item The commutativity of the face labeled $\Cc$ (notice that this face
fully encloses the one labeled $\Ce$) is due to
Duflo~\cite[th\'eor\`eme~V.1]{Duflo:Resolubles}.

\end{myitemize}

\begin{proposition} \lbl{PreciseWheelingTheorem}
The face labeled $\Cd$ is commutative.
\end{proposition}

\begin{remark} Recalling that $\Dj$ is an algebra isomorphism,
Proposition~\ref{PreciseWheelingTheorem} becomes the precise formulation
of Theorem~\ref{WheelingTheorem}.
\end{remark}

\begin{proof}[Proof of Proposition~\ref{PreciseWheelingTheorem}] Follows
immediately from the following two lemmas, taking $C=\Omega$ in
\eqref{ChatC}.
\def\qed{$\hfill\smily$}
\end{proof}

\begin{lemma} \lbl{FirstLemma} Let $\kappa:\frakg\to\frakg^\star$ be the
identification induced by the standard bilinear form $(\cdot,\cdot)$ of
$\frakg$. Extend $\kappa$ to all symmetric powers of $\frakg$, and let
$\kappa^\hbar:S(\frakg)^\frakg[[\hbar]]\to S(\frakg^\star)(\!(\hbar)\!)$
be defined for a homogeneous $s\in S(\frakg)^\frakg[[\hbar]]$ (relative to
the grading of $S(\frakg)$) by $\kappa^\hbar(s)=\hbar^{-\deg s}\kappa(s)$.
If $C\in\calB'$ is a Chinese character, $\hat{C}:\calB'\to\calB'$ is
the operator corresponding to $C$ as in Definition~\ref{ChatDef}, and
$C'\in\calB'$ is another Chinese character, then
\begin{equation} \lbl{ChatC}
  \Tg\hat{C}(C') = D(\kappa^\hbar\Tg C)\Tg C'.
\end{equation}
\end{lemma}

\begin{proof} If $\kappa j$ is a tensor in
$S^{k}(\frakg^\star)\subset\frakg^{\star\otimes k}$, the $k$'th
symmetric tensor power of $\frakg^\star$, and $j'$ is a tensor in
$S^{k'}(\frakg)\subset\frakg^{\otimes k'}$, then
\begin{equation} \lbl{DInFull}
  D(\kappa j)(j')=\begin{cases}
    0 & \text{if $k>k'$,} \\
    \parbox{2.6in}{
      the sum of all ways of contracting all the tensor components of $j$
      with some (or all) tensor components of $j'$
    }\quad & \text{otherwise.}
  \end{cases}
\end{equation}
By definition, the ``diagrams to Lie algebras'' map carries gluing to
contraction, and hence carries the operation in Definition~\ref{ChatDef}
to the operation in~\eqref{DInFull}, namely, to $D$. Counting powers of
$\hbar$, this proves~\eqref{ChatC}.
\def\qed{$\hfill\smily$}
\end{proof}

\begin{lemma} \lbl{TgOmega}
$\kappa^\hbar\Tg\Omega=j_\frakg^{1/2}$.
\end{lemma}

\begin{proof} It follows easily from the definition of $\Tg$ and $\kappa^h$
that $(\kappa^\hbar\Tg\omega_n)(X)=\tr(\ad X)^n$ for any $X\in\frakg$.
Hence, using the fact that $\kappa^\hbar\circ\Tg$ is an algebra morphism if
$\calB'$ is taken with the disjoint union product,
\[ (\kappa^\hbar\Tg\Omega)(X)
  = \exp\sum_{n=1}^\infty b_{2n}(\kappa^\hbar\Tg\omega_{2n})(X)
  = \exp\sum_{n=1}^\infty b_{2n}\tr(\ad X)^{2n}
  = \det\exp\sum_{n=1}^\infty b_{2n}(\ad X)^{2n}.
\]
By the definition of the modified Bernoulli numbers~\eqref{MBNDefinition},
this is
\begin{equation}
  \det\exp\frac{1}{2}\log\frac{\sinh \ad X/2}{\ad X/2}
  = \det\left(\frac{\sinh \ad X/2}{\ad X/2}\right)^{1/2}
  = j_\frakg^{1/2}(X).
  \prtag{\smily}
\end{equation}
\renewcommand{\qed}{}
\end{proof}

\begin{proposition} 
The face labeled $\Ce$ is commutative.
\end{proposition}

\begin{proof} According to M.~Vergne (private communication), this is a
well known fact. We could not find a reference, so here's the gist of
the proof. Forgetting about powers of $\hbar$ and $\frakg$-invariance
and taking the Fourier transform (over $\frakg_{\mathbf R}$), the
differential operator $\Dj$ becomes the operator of multiplication by
$j_\frakg^{1/2}(iX)$ on $S(\frakg)$. Taking the inverse Fourier transform,
we see that $\Dj$ is the operator of convolution with the inverse Fourier
transform of $j_\frakg^{1/2}(iX)$, which is the symplectic measure on
$M_\rho$ (see~\eqref{FourierTransform}). So $\Dj$ is convolution with
that measure, as required.
\def\qed{$\hfill\smily$}
\end{proof}

\section{Proof of Theorem~\ref{UnknotTheorem}} \lbl{proof}

We prove the slightly stronger equality
\begin{equation} \lbl{PreciseUnknotEquation}
  {\mathcal T}^\hbar_\frakg\Omega
  = {\mathcal T}^\hbar_\frakg\chi^{-1}\ZF(\bigcirc).
\end{equation}

\begin{proof} We compute the right hand side
of~\eqref{PreciseUnknotEquation} by computing
$\Sg\ig^{-1}\psig\RTg(\bigcirc)$ and using the commutativity of the monster
diagram. It is known (see
e.g.~\cite[example~11.3.10]{ChariPressley:QuantumGroups}) that if
$\lambda-\rho\in\frakh^\star$ is the highest weight of some irreducible
representation $R_{\lambda-\rho}$ of $\frakg$, then
\[ (\psig\RTg(\bigcirc))(\lambda)
  = \frac{1}{\dim R_{\lambda-\rho}} \prod_{\alpha\in\Delta_+} \frac
    {\sinh\hbar(\lambda,\alpha)/2}
    {\sinh\hbar(\rho,\alpha)/2},
\]
where $\Delta_+$ is the set of positive roots of $\frakg$ and
$(\cdot,\cdot)$ is the standard invariant bilinear form on $\frakg$.
By the Weyl dimension formula and some minor arithmetic, we get (see
also~\cite[section~7]{LeMurakami:Parallel})
\begin{equation} \lbl{ProductOverRoots}
  (\psig\RTg(\bigcirc))(\lambda) = \prod_{\alpha\in\Delta_+}
  \frac
    {\hbar(\rho,\alpha)/2}
    {\sinh\hbar(\rho,\alpha)/2}
  \cdot\frac
    {\sinh\hbar(\lambda,\alpha)/2}
    {\hbar(\lambda,\alpha)/2}
  .
\end{equation}
We can identify $\frakg$ and $\frakg^\star$ using the form $(\cdot,\cdot)$,
and then expressions like `$\ad\lambda$' makes sense. By definition,
if $\frakg_\alpha$ is the weight space of the root $\alpha$, then
$\ad\lambda$ acts as multiplication by $(\lambda,\alpha)$
on $\frakg_\alpha$, while acting trivially on $\frakh$. From this
and~\eqref{ProductOverRoots} we get
\[ (\psig\RTg(\bigcirc))(\lambda) =
  \det\left(\frac{\ad\hbar\rho/2}{\sinh\ad\hbar\rho/2}\right)^{1/2}
  \cdot
  \det\left(\frac{\sinh\ad\hbar\lambda/2}{\ad\hbar\lambda/2}\right)^{1/2} =
  j_\frakg^{-1/2}(\hbar\rho)\cdot j_\frakg^{1/2}(\hbar\lambda).
\]

The above expression (call it $Z(\lambda)$) makes sense
for all $\lambda\in\frakg^\star$, and hence it is also
$\ig^{-1}\psig\RTg(\bigcirc)$. So we're only left with computing $\Sg
Z(\lambda)$:
\[ \Sg Z(\lambda)=\int_{r\in M_{i\rho}}dr\,Z(\lambda-ir)=
  j_\frakg^{-1/2}(\hbar\rho)
    \int_{r\in M_{i\rho}}dr\,j_\frakg^{1/2}(\hbar(\lambda-ir)).
\]
By~\eqref{FourierTransform}, this is
\[ j_\frakg^{-1/2}(\hbar\rho)\int_{r\in M_{i\rho}}dr\int_{r'\in M_{i\rho}}dr'\,
  e^{i\hbar( r',\lambda-ir)}
  = j_\frakg^{-1/2}(\hbar\rho)\int_{r'\in M_{i\rho}}dr'\,
    e^{i\hbar( r',\lambda)}\int_{r\in M_{i\rho}}dr\,
    e^{i\hbar( -ir',r)}.
\]
Using~\eqref{FourierTransform} again, we find that the inner-most integral
is equal to $j_\frakg^{1/2}(\hbar\rho)$ independently of $r'$, and hence
\[ \Sg Z(\lambda)
  = \int_{r'\in M_{i\rho}}dr'\,e^{i\hbar( r',\lambda)},
\]
and using~\eqref{FourierTransform} one last time we find that
\begin{equation} \lbl{SgZ}
  \Sg Z(\lambda) = j_\frakg^{1/2}(\hbar\lambda).
\end{equation}

The left hand side of~\eqref{PreciseUnknotEquation} was already computed
(up to duality and powers of $\hbar$) in Lemma~\ref{TgOmega}. Undoing the
effect of $\kappa^\hbar$ there, we get the same answer as in~\eqref{SgZ}.
\def\qed{$\hfill\smily$}
\end{proof}

\[ \printname{WheelsWheeling}
	\setlength{\unitlength}{0.25\standardunitlength}
	\begin{array}{c}  \hspace{-1.7mm}
        	\raisebox{-8pt}{\input draws/WheelsWheeling.tex }
        	\hspace{-1.9mm}
	\end{array}
 \]

\ifx\undefined\bysame
        \newcommand{\bysame}{\leavevmode\hbox to3em{\hrulefill}\,}
\fi


\end{document}